\documentclass[11pt,a4paper]{article}
\usepackage[T1]{fontenc}
\usepackage{lmodern,amsmath,amssymb,amsthm,microtype}
\usepackage[margin=26mm]{geometry}
\usepackage[hidelinks]{hyperref}
\newtheorem{theorem}{Theorem}
\newtheorem{lemma}[theorem]{Lemma}
\newtheorem{proposition}[theorem]{Proposition}
\newtheorem{corollary}[theorem]{Corollary}
\theoremstyle{definition}
\newtheorem{definition}[theorem]{Definition}
\newtheorem{problem}[theorem]{Problem}
\DeclareMathOperator{\HG}{HG}
\newcommand{\Pp}{\mathbb P}
\newcommand{\Ee}{\mathbb E}

\title{Sublinear Bounds for Hat Guessing Numbers on Bipartite Graphs and Directed Blow-ups}
\author{\v{Z}arko Ran\dj elovi\'c*}
\date{November 2026}
\begin{document}
\maketitle
\begin{abstract}
Consider a game in which players are placed at the vertices
of a graph, each wearing a hat of one of $q$ possible
colours. Each player can see the hats of their neighbours
but not their own, and the players simultaneously guess
the colours of their own hats. The hat-guessing number
of the graph is the largest $q$ for which the players
have a strategy that guarantees at least one correct
guess, regardless of the assignment of hats.
Gadouleau and Georgiou showed that the hat-guessing
number of the complete bipartite graph $K_{n,n}$ is
at most $n+1$. We show that it is in fact at most
$(1+o(1))n/\log\log n$, hence proving a conjecture
of He, Ido and Przybocki for sufficiently large $n$. More generally, we show that the hat-guessing number 
$\mathrm{HG}(K_{m,n})=o(m)$ whenever $m\leq n$ and
$\log n=o(m)$, as $m\to\infty$. We extend our results
to certain classes of directed graphs, including
blow-ups of directed cycles of arbitrary length.
We also show that the leading asymptotic term in our
bound cannot be improved within the probabilistic
framework used in this paper.
\end{abstract}

\section{Introduction}
Suppose that $n$ players are each given a hat whose colour is chosen from a set of $q$ colours. Each player can see the hats of some of the other players, but cannot see their own hat. Before the hats are assigned, the players agree on a strategy; once they have seen the hats, they must simultaneously guess the colour of their own hat. They win if at least one player guesses correctly. How large can $q$ be if the players want to guarantee a win, regardless of the assignment of hats?

The answer clearly depends on which hats each player can see. We describe this by a finite loopless directed graph $G$, where an arc $v\to w$ means that the player at $v$ sees the hat of the player at $w$. The largest number of colours for which a winning strategy exists is called the \emph{hat-guessing number} of $G$, denoted by $\HG(G)$. In the case of undirected graphs, each edge corresponds to two arcs of opposite orientations. This version of the hat-guessing game was introduced by Butler, Hajiaghayi, Kleinberg and Leighton \cite{butler}.

Let $Q$ be the set of available colours, and for a finite
set $U$, write $Q^U$ for the set of all colour
assignments on $U$. For each
vertex $v$ of $G$, let $N^+(v)$ denote the set of
out-neighbours of $v$. A \emph{guessing function}
for the player at $v$ is a function
$f_v:Q^{N^+(v)}\to Q$, where the
player guesses $f_v(x)$ when they see the colour assignment $x$.
A strategy consists of a choice of a guessing
function for each vertex of $G$.

\makeatletter{\let\thefootnote\relax\footnote{*Mathematical Institute of the Serbian Academy of Sciences and Arts, Kneza Mihaila 36, Belgrade 11000, Serbia. Email: zarko.randjelovic@turing.mi.sanu.ac.rs}}

If every player can see all the other hats, then $\HG(K_n)=n$. Indeed, when the colours are identified with the residues modulo $n$, one can have the $i$-th player guess the colour that would make the sum of all the colours equal to $i$ modulo $n$. Whatever the assignment, exactly one of the guesses is correct. On the other hand, for $q>n$, if the hats are assigned independently and uniformly, the probability that a particular player guesses correctly is $1/q$, so the probability that at least one player guesses correctly is at most $n/q<1$.

The situation is rather different when the players are split into two groups and each player can see precisely the hats of the players in the other group. This is the game on the complete bipartite graph. Initially, Butler, Hajiaghayi, Kleinberg and Leighton
\cite{butler} obtained a lower bound of order
$\frac{\log\log n}{\log\log\log n}$ for $\HG(K_{n,n})$. Gadouleau and Georgiou \cite{gadouleau} improved
this to a lower bound of order $\frac{\log n}{\log\log n}$
and showed that
$\HG(K_{m,n})\leq\min\{m,n\}+1$. More recently, Alon, Ben-Eliezer, Shangguan and Tamo \cite{alon} substantially improved the lower bound in the balanced case, so that the two bounds are
\[
 n^{1/2-o(1)}\le \HG(K_{n,n})\le n+1.
\]
He, Ido and Przybocki \cite{he} showed that $\HG(K_{3,3})=3$ and conjectured that $\HG(K_{n,n})\le n$ for every $n\ge3$. They also suggested that the hat-guessing number of $K_{n,n}$ might be sublinear. In this paper, we show that this is indeed the case. General upper bounds in terms of graph order and maximum degree have recently been studied by Shurman and Sibley \cite{shurman}.

The upper bound of $\min\{m,n\}+1$ for $\HG(K_{m,n})$
is actually attained for every fixed $m$ when $n$ is
sufficiently large (see \cite{gadouleau}). However, the situation is quite
different when the sizes of the two parts do not
differ too much.
We show that $\HG(K_{m,n})=o(m)$ as $m\to\infty$
whenever $m\leq n$ and $\log n=o(m)$.
More precisely, we prove the following result.

\begin{theorem}\label{main}
Let $m\le n$, with $m\to\infty$ and $\log n=o(m)$. Then
\[
 \HG(K_{m,n})\le (1+o(1))\frac{m}{\log\log(m/\log n)}.
\]
More precisely, writing $A=\log\log(m/\log(2n))$, we have, for sufficiently large $m$,
\begin{equation}\label{fine}
 \HG(K_{m,n})<\frac{m}{A-\log A-\log\log A-3}.
\end{equation}
\end{theorem}
Here and throughout the paper, all logarithms are natural. In particular, Theorem~\ref{main} gives
\[
 \HG(K_{n,n})\le(1+o(1))\frac{n}{\log\log n}=o(n).
\]
The theorem applies even when the two parts have very different sizes. For example, if $n=\lfloor e^{\sqrt m}\rfloor$, then the upper bound is asymptotically the same as for $K_{m,m}$. Since adding edges cannot decrease the hat-guessing number, we also obtain a sublinear upper bound for all bipartite graphs in terms of the total number of vertices.

We briefly describe the idea of the proof. Suppose we are given the guessing functions of the players in one part of $K_{m,n}$. We want to find a set of assignments of hats in the other part with two properties: no guessing function takes all $q$ values on this set, and, regardless of what the players in the other part guess, there is an assignment from the set on which all their guesses are incorrect. The difficulty is that the second property favours large sets, whereas the first favours small sets. By restricting the possible assignments in a suitable way before choosing a random set, we show that both properties can be satisfied simultaneously.

The same construction turns out to be useful for certain directed graphs. In particular, if each vertex of a directed cycle is replaced by $n$ independent vertices, its hat-guessing number is at most $(1+o(1))n/\log\log n$, regardless of the length of the cycle.

The paper is organised as follows.
In Section~\ref{sec:sampling} we prove
Theorem~\ref{main}.
In Section~\ref{sec:directed} we consider directed
graphs and give a more general sufficient condition
for obtaining sublinear upper bounds.
Finally, in Section~\ref{sec:barrier} we show that our asymptotic bound on
$\HG(K_{m,n})$ cannot be improved within the
probabilistic framework used in the proof of
Theorem~\ref{main}.

\section{Proof of Theorem~\ref{main}}\label{sec:sampling}
Let $Q$ be a set of $q$ colours, and write
$[k]=\{1,\ldots,k\}$. For $z\in Q^U$, define the cube with centre $z$ by
\[
 C_z=\{x\in Q^U:x_u\ne z_u\text{ for every }u\in U\}.
\]
A subset $S\subseteq Q^U$ is \emph{blocking} if it meets every $C_z$.

We will use two simple observations. First, note that adding arcs to $G$ cannot decrease $\HG(G)$.
Secondly, to show that $\HG(G)<q$, it suffices to show
that there is no winning strategy with $q$ colours (since reducing the number of colours can only help
the players). 

We first establish a connection between blocking
sets and winning strategies. The following lemma
is a reformulation of the Hamming ball condition
used in \cite[Lemma~3.1]{alon}, which goes back
to \cite{gadouleau}; see also
\cite[Propositions~5.1.2 and~5.1.3]{mcinnis}..

\begin{lemma}\label{crit}
Let $q\ge2$, and fix the $n$ right-side guessing functions $f_i:Q^m\to Q$ of $K_{m,n}$. They cannot be completed to a winning strategy if and only if there exists a blocking set $S\subseteq Q^m$ such that $f_i(S)\ne Q$ for all $i$.
\end{lemma}
\begin{proof}
Suppose first that such an $S$ exists. For each $i$, choose $y_i\in Q\setminus f_i(S)$ and give the $i$-th player on the right a hat of colour $y_i$. If the assignment on the left belongs to $S$, all the players on the right guess incorrectly. Now the guesses of the players on the left are determined by $y=(y_i)_{i=1}^n$; denote their vector of guesses by $z$. Since $S$ is blocking, there is an $x\in S\cap C_z$. With this assignment on the left, all the players guess incorrectly.

Conversely, suppose that there is no blocking set
$S\subseteq Q^m$ such that $f_i(S)\ne Q$ for every
$i\in[n]$. For each $y\in Q^n$, define $S_y=\{x:f_i(x)\ne y_i\ \text{ for all } i\}$. None of these sets are blocking, so for each $y$ choose $z(y)$ with $S_y\cap C_{z(y)}=\varnothing$. The players on the left can use $z(y)$ as their vector of guesses when the colour assignment on
the right is $y$. Since $S_y\cap C_{z(y)}=\varnothing$, at least one guess is correct whenever all right guesses are wrong. This completes a winning strategy.
\end{proof}

As an immediate consequence of Lemma~\ref{crit},
we recover the upper bound of Gadouleau and Georgiou
\cite{gadouleau}. Indeed, if $q\geq m+2$, then any
set of $m+1$ constant vectors of distinct colours
in $Q^m$ is blocking. Since this set has size less
than $q$, Lemma~\ref{crit} gives
$\HG(K_{m,n})\leq m+1$.
Exchanging the roles of the two parts, we obtain
$\HG(K_{m,n})\leq\min\{m,n\}+1$.

We now turn to the construction of blocking sets. We will first prove a statement about an arbitrary family of functions on $Q^d$.
Let $d$ be a positive integer and $q\ge2$ an integer. Partition the $d$ coordinates into
$g=\lceil d/k\rceil$ nonempty blocks
$B_1,\ldots,B_g$, each of size at most $k$. Let $\Omega\subseteq Q^d$ consist of all assignments constant on each block, and put $N=|\Omega|=q^g$. For a cube center $z$, the forbidden set in group $j$ is $A_j=\{z_i:i\in G_j\}$. Then
\begin{equation}\label{cubeproduct}
 \frac{|C_z\cap\Omega|}{N}=\prod_{j=1}^g(1-|A_j|/q),\qquad \sum_j|A_j|\le d.
\end{equation}
Notice that $C_z\cap\Omega$ is determined entirely
by the sets $A_1,\ldots,A_g$. We will refer to $C_z\cap\Omega$
the restriction of the cube $C_z$ to $\Omega$ and
define its density to be $|C_z\cap\Omega|/|\Omega|$. The following lemma bounds the number
of distinct cube restrictions and gives a lower
bound on their sizes.

\begin{lemma}\label{count}
Suppose that $2\leq q\leq d$, put $t=d/q$, and assume
that $t\geq 2$ and $q/t\geq 2$. Let $\Omega$ be
constructed as above with $k=\lfloor q/t\rfloor$.
Then
\[
    \frac{|C_z\cap\Omega|}{|\Omega|}
    \geq e^{-t-2}
\]
for every $z\in Q^d$, and there are at most
$(3et)^d$ distinct cube restrictions to $\Omega$.
\end{lemma}
\begin{proof}
The inequality $-\log(1-u)\leq u/(1-u)$ holds for
$0\leq u<1$. Indeed, setting $x=1-u$, this is
equivalent to
\[
    x\log x+1-x\geq 0
\]
for $0<x\leq 1$, which follows by differentiation. Since $|A_j|/q\leq 1/t$ for every $j\in[g]$,
it follows from \eqref{cubeproduct} that
\[
 -\log(|C_z\cap\Omega|/N)\le \frac{1}{1-1/t}\sum_j |A_j|/q
 \le \frac{t^2}{t-1}\le t+2.
\]
This proves the first assertion of the lemma. It remains to bound the number of distinct restrictions. Each $g$-tuple $(A_j)_{j=1}^g$ corresponds to a subset of $[g]\times Q$ of size at most $d$. There are at most $\sum_{h=0}^d\binom{gq}{h}$ such subsets. Since $gq\ge d$, take $u=d/(gq)\le1$. The binomial theorem gives
\begin{equation}\label{binomialbound}
    u^d\sum_{h=0}^d\binom{gq}{h}
    \leq (1+u)^{gq}\leq e^{ugq}=e^d.
\end{equation}
Since $q\geq 2t$, we have
$k=\lfloor q/t\rfloor\geq q/t-1\geq q/(2t)$.
Hence,
\[
    gq\leq\left(\frac{d}{k}+1\right)q
    \leq 2dt+q\leq 3dt.
\] Dividing both sides of \eqref{binomialbound} by $u^d$,
we obtain
\[
    \sum_{h=0}^d\binom{gq}{h}
    \leq \left(\frac{e}{u}\right)^d
    = \left(\frac{egq}{d}\right)^d
    \leq (3et)^d,
\]
as required.
\end{proof}

Having bounded both the number and the densities of the cube restrictions, we now choose a random subset
of $\Omega$, constructed as in Lemma~\ref{count}. Let $M>0$ be a parameter to be chosen later.
Select each point of $\Omega$ independently with
probability $\theta=1-e^{-M/N}$, forming a random
set $S$. For any fixed $B\subseteq\Omega$, we have
\begin{equation}\label{avoid}
    \Pp(S\cap B=\varnothing)=\exp(-M|B|/N).
\end{equation}
The point of choosing this probability is that
the probability of missing a set $B$ now has a
particularly simple expression. Next, we estimate
the probability that a fixed function
$f:Q^d\to Q$ takes every value in $Q$ on the
random set $S$.

\begin{lemma}\label{surj}
For any fixed $f:Q^d\to Q$ and any $M>0$, the above random set satisfies
\begin{equation}\label{surjformula}
 \Pp(f(S)=Q)\le \exp(-q e^{-M/q}).
\end{equation}
\end{lemma}
\begin{proof}
Let $B_a=\{x\in\Omega:f(x)=a\}$ and $\mu_a=|B_a|/N$. The sets $B_a$ partition $\Omega$, and $\sum_a\mu_a=1$. The events that $S$ meets $B_a$, for $a\in Q$,
are mutually independent, since the sets $B_a$
are pairwise disjoint and each point of $\Omega$
is selected independently. By \eqref{avoid} and and using the elementary inequality $1-u\leq e^{-u}$, we obtain
\[
 \Pp(f(S)=Q)=\prod_{a\in Q}(1-e^{-M\mu_a})\le
 \exp\left(-\sum_{a\in Q}e^{-M\mu_a}\right).
\]
The AM-GM inequality gives $\sum_a e^{-M\mu_a}\ge q(\prod_a e^{-M\mu_a})^{1/q}=q e^{-M/q}$. This proves \eqref{surjformula}. In particular, the upper bound does not depend on the function $f$.
\end{proof}

\begin{theorem}\label{family}
Suppose the hypotheses of Lemma~\ref{count} hold. Define
\begin{equation}\label{M}
 M_d=e^{t+2}\bigl(d\log(3et)+\log4\bigr).
\end{equation}
For any family of $r\ge1$ functions $f_1,\ldots,f_r:Q^d\to Q$, if
\begin{equation}\label{test}
 q e^{-M_d/q}\ge \log(2r),
\end{equation}
there exists a blocking set $S$ such that every $f_i(S)$ is a proper subset of $Q$. 
\end{theorem}
\begin{proof}
Let $M=M_d$. We first estimate the probability that $S$ is not blocking. By Lemma~\ref{count} and \eqref{avoid}, summing over the distinct restrictions $C_z\cap\Omega$, we obtain
by the union bound
\[
 \Pp(S\text{ is not blocking})\le (3et)^d\exp(-M_de^{-t-2})=1/4.
\]
By Lemma~\ref{surj}, \eqref{test} and the union bound, the probability that at least one of the $r$ functions is surjective on $S$ is at most $r\exp(-qe^{-M_d/q})\le1/2$. Consequently, with probability at least $1/4$, the set $S$ is blocking and none of the functions $f_i$ take all values in $Q$ on $S$. In particular, such an $S$ exists.
\end{proof}

We will also use this result for functions on fewer than $d$ coordinates. Indeed, we may extend the functions to $Q^d$ by having them ignore the extra coordinates and then project the blocking set obtained from the theorem onto the original coordinates. Projection preserves blocking, because every cube center on the smaller coordinate set extends to a cube center on the full set. Since none of the extended functions take all values
in $Q$ on the blocking set, none of the original functions do so
on the projected set.

Before proving Theorem~\ref{main}, we compare the bounds one obtains by constructing blocking sets on the two different sides of a complete bipartite graph.
We may apply Theorem~5 with either $d=m$ and $r=n$,
or $d=n$ and $r=m$, and then use Lemma~\ref{crit}
to obtain sufficient conditions for $\HG(K_{m,n})<q$. To facilitate the comparison between these two
approaches, we give a simpler, though less sharp,
sufficient condition that does not require the
additional assumptions $t\geq 2$ and $q/t\geq 2$.

\begin{proposition}\label{orientation}
For $2\leq q\leq d$, put $
    M_d^{(0)}=(4d\log 2+\log 4)e^{2d/q}.$
For any family of $r\geq 1$ functions
$f_1,\ldots,f_r:Q^d\to Q$, the conclusion of
Theorem~\ref{family} holds if
\begin{equation}\label{test0}
    qe^{-M_d^{(0)}/q}\geq \log(2r).
\end{equation}
Moreover, if $q\leq m\leq n$ and \eqref{test0}
holds for $(d,r)=(n,m)$, then it also holds for
$(d,r)=(m,n)$.
\end{proposition}

\begin{proof}
We use the same construction of $\Omega$ as before,
but with $k=\lfloor q/2\rfloor$. Since
$|A_j|/q\leq 1/2$ for all $j\in[g]$, it follows
from \eqref{cubeproduct} and the inequality
$-\log(1-u)\leq 2u$ for $0\leq u\leq 1/2$ that
\[
    -\log\frac{|C_z\cap\Omega|}{|\Omega|}
    \leq 2\sum_{j=1}^g\frac{|A_j|}{q}
    \leq \frac{2d}{q}.
\]
Hence, for all $z\in Q^d$,
\[
    \frac{|C_z\cap\Omega|}{|\Omega|}
    \geq e^{-2d/q}.
\]

Since $k=\lfloor q/2\rfloor\geq q/3$, we have
\[
    gq\leq\left(\frac{d}{k}+1\right)q
    \leq 3d+q\leq 4d.
\]
Recall that each $g$-tuple $(A_1,\ldots,A_g)$
corresponds to a subset of $[g]\times Q$.
Thus, the number of distinct cube restrictions
to $\Omega$ is at most $2^{gq}\leq 2^{4d}$.

We now choose a random subset $S$ of $\Omega$
as before, with $M=M_d^{(0)}$.
By \eqref{avoid} and the union bound, we obtain
\[
    \Pp(S\text{ is not blocking})
    \leq 2^{4d}\exp(-M_d^{(0)}e^{-2d/q})
    =\frac14.
\]
On the other hand, by Lemma~\ref{surj},
\eqref{test0} and the union bound, we have
\[
    \Pp(f_i(S)=Q\text{ for at least one  }i\in[r])
    \leq r\exp(-qe^{-M_d^{(0)}/q})
    \leq \frac12.
\]
As in the proof of Theorem~\ref{family}, with
positive probability $S$ is blocking and none
of the functions $f_i$ take all values in $Q$
on $S$.

It remains to prove the second assertion.
Suppose that $q\leq m\leq n$ and that
\[
    qe^{-M_n^{(0)}/q}\geq\log(2m).
\]
By the definition of $M_d^{(0)}$, we have
\begin{align*}
    M_n^{(0)}-M_m^{(0)}
    &\geq
    (4n\log 2+\log 4)e^{2m/q}
    -(4m\log 2+\log 4)e^{2m/q}\\
    &=4(n-m)\log 2\cdot e^{2m/q}\\
    &\geq n-m.
\end{align*}
Consequently, since $q\leq m$, we obtain
\[
    \exp\left(\frac{M_n^{(0)}-M_m^{(0)}}{q}\right)
    \geq e^{(n-m)/m}\geq\frac{n}{m},
\]
where the last inequality follows from
$e^x\geq 1+x$.

It follows that
\begin{align*}
    qe^{-M_m^{(0)}/q}
    &=qe^{-M_n^{(0)}/q}
    \exp\left(\frac{M_n^{(0)}-M_m^{(0)}}{q}\right)\\
    &\geq\frac{n}{m}\log(2m)\\
    &\geq\log(2n).
\end{align*}
Indeed, the last inequality follows from
$\log(n/m)\leq n/m-1$ and $\log(2m)\geq 1$.
Thus, \eqref{test0} also holds for
$(d,r)=(m,n)$, as required.
\end{proof}

Proposition~\ref{orientation} shows that, for the
sufficient condition given in \eqref{test0}, there
is no advantage in constructing a blocking set
on the larger part of a complete bipartite graph.
We are now ready to prove Theorem~\ref{main}.

\begin{proof}[Proof of Theorem~\ref{main}]
Let $A$ be as in the statement of the theorem,
and put $R=e^A$. In particular,
\[
    R=\log(m/\log(2n)).
\]
Since $\log n=o(m)$, we have $R\to\infty$ and
hence $A\to\infty$ as $m\to\infty$.
Moreover, $A\leq\log\log m$ for sufficiently
large $m$.

Set
\[
    T=A-\log A-\log\log A-3,
    \qquad
    q=\left\lceil\frac{m}{T}\right\rceil,
    \qquad
    t=\frac{m}{q}.
\]
We will show that there is no winning strategy
with $q$ colours.

Since $A\to\infty$, we have $T\sim A$.
Also, since $A\leq\log\log m$, we obtain
$q\sim m/A$ and $t\sim A$.
By the choice of $q$, we have $t\leq T$.
In particular, for sufficiently large $m$,
the hypotheses of Lemma~\ref{count} hold with
$d=m$.

We next verify that \eqref{test} holds with
$d=m$ and $r=n$. By \eqref{M}, we have
\begin{align*}
    \frac{M_m}{q}
    &=e^{t+2}t
    \left(\log(3et)+\frac{\log 4}{m}\right)\\
    &\leq e^{T+2}T
    \left(\log(3eT)+\frac{\log 4}{m}\right).
\end{align*}
By the definition of $T$,
\[
    e^{T+2}
    =\frac{e^{-1}e^A}{A\log A}
    =\frac{e^{-1}R}{A\log A}.
\]
Therefore,
\[
    \frac{M_m}{q}
    \leq e^{-1}R\cdot\frac{T}{A}\cdot
    \frac{\log(3eT)+m^{-1}\log 4}{\log A}.
\]
Since $T\sim A$, both fractions on the
right-hand side tend to $1$. Consequently,
\begin{equation}\label{Mmestimate}
    \frac{M_m}{q}\leq(e^{-1}+o(1))R.
\end{equation}

On the other hand, by the choice of $q$,
\begin{align*}
    \log\left(\frac{q}{\log(2n)}\right)
    &\geq
    \log\left(\frac{m}{T\log(2n)}\right)\\
    &=R-\log T\\
    &=(1+o(1))R,
\end{align*}
where the last equality follows from
$T\sim A=\log R$.

It follows from
\eqref{Mmestimate} that, for sufficiently large
$m$,
\[
    \frac{M_m}{q}
    \leq\log\left(\frac{q}{\log(2n)}\right).
\]
Equivalently,
\[
    qe^{-M_m/q}\geq\log(2n).
\]
Thus, \eqref{test} holds with $d=m$ and $r=n$.

By Theorem~\ref{family}, for any choice of the
$n$ right-side guessing functions
$f_1,\ldots,f_n:Q^m\to Q$, there exists a
blocking set $S\subseteq Q^m$ such that
$f_i(S)\ne Q$ for all $i\in[n]$.
By Lemma~\ref{crit}, there is no winning strategy
with $q$ colours. Hence,
\[
    \HG(K_{m,n})
    \leq q-1
    <\frac{m}{T}
    =\frac{m}{A-\log A-\log\log A-3},
\]
proving \eqref{fine}.

Finally, since $T\sim A$ and $
    \log\log(m/\log(2n))
    \sim\log\log(m/\log n)$
we obtain
\[
    \HG(K_{m,n})
    \leq(1+o(1))
    \frac{m}{\log\log(m/\log n)},
\]
as required.
\end{proof}

As a consequence of Theorem~\ref{main}, we obtain
a simpler asymptotic bound when the size of the
larger part is not too large.

\begin{corollary}\label{balanced}
As $n\to\infty$, we have
\[
    \HG(K_{n,n})\leq
    (1+o(1))\frac{n}{\log\log n}.
\]
More generally, if $m\leq n$ and
$\log n=O(m^{1-\delta})$ for some fixed
$\delta>0$, then
\[
    \HG(K_{m,n})\leq
    (1+o(1))\frac{m}{\log\log m}
\]
as $m\to\infty$.
\end{corollary}

\begin{proof}
Since $\log n=O(m^{1-\delta})$, we have
$\log n=o(m)$ and
\[
    \delta\log m-O(1)
    \leq \log(m/\log n)
    \leq \log m
\]
for sufficiently large $m$. Hence,
$
    \log\log(m/\log n)=\log\log m+O(1).
$
The second assertion now follows from
Theorem~\ref{main}, and the first follows
by taking $m=n$.
\end{proof}

For example, Corollary~\ref{balanced} applies
when $n=m^{10}$ or when
$n=\lfloor e^{\sqrt m}\rfloor$.
If instead $n=\lfloor e^{m/\log m}\rfloor$,
Theorem~\ref{main} gives
\[
    \HG(K_{m,n})
    \leq(1+o(1))\frac{m}{\log\log\log m}.
\]

Since adding edges cannot decrease the hat-guessing
number, Theorem~\ref{main} also gives a sublinear
upper bound for all bipartite graphs.

\begin{corollary}\label{allbip}
Every bipartite graph $G$ on $N$ vertices satisfies
\[
    \HG(G)\leq(1+o(1))\frac{N}{2\log\log N}
\]
as $N\to\infty$, where the error term is independent
of $G$. The same bound holds for directed graphs
whose underlying undirected graphs are bipartite.
\end{corollary}

\begin{proof}
Let $m$ and $n$ be the sizes of the two parts
of $G$, where $m\leq n$ and $m+n=N$.
The case $m=0$ is trivial so we may assume $m\ge 1$. We know that, $\HG(G)\leq\HG(K_{m,n})$.

If $m\leq\sqrt N$, the bound
$\HG(K_{m,n})\leq m+1$ gives
\[
    \HG(G)\leq\sqrt N+1
    =o\left(\frac{N}{\log\log N}\right).
\]
We may therefore assume that $m>\sqrt N$.
In this case, $\log n\leq\log N=o(m)$.
Furthermore,
\[
    \frac{\sqrt N}{\log N}
    <\frac{m}{\log n}\leq N
\]
for sufficiently large $N$. It follows that
\[
    \log\log(m/\log n)
    =\log\log N+O(1).
\]
By Theorem~\ref{main} and the fact that
$m\leq N/2$, we obtain
\[
    \HG(G)\leq\HG(K_{m,n})
    \leq(1+o(1))\frac{m}{\log\log N}
    \leq(1+o(1))\frac{N}{2\log\log N}.
\]

The directed case follows immediately by
adding all possible arcs in both directions.
\end{proof}

\section{Directed graphs and blow-ups}
\label{sec:directed}

In this section, we extend our results to certain
classes of directed graphs.

Recall that a directed graph is acyclic if it
contains no directed cycle.
The vertices of an acyclic directed graph can
be ordered so that every arc goes from an
earlier vertex to a later vertex. Such an
ordering is called a topological ordering.

It is well known that every nonempty acyclic
directed graph has hat-guessing number $1$
(see, for example, \cite[Corollary~2.2.7]{mcinnis}).
Indeed, suppose that $q\geq 2$ and fix an
arbitrary strategy.
Assign hats to the vertices in reverse
topological order.
When assigning a hat to a vertex, the colours
of all the hats it can see have already been
chosen, so its guess is determined.
We can therefore choose its hat to have a
different colour from its guess.
In this way, we obtain an assignment on which
all the players guess incorrectly.

We will use the following known result, which
shows that the hat-guessing number of a directed
graph is determined by the hat-guessing numbers
of its strongly connected components.

\begin{proposition}[\cite{mcinnis}, Proposition~2.2.6]
\label{scc}
For every nonempty directed graph $G$, we have
\[
    \HG(G)=\max_C\HG(G[C]),
\]
where the maximum is taken over all strongly
connected components $C$ of $G$, and $G[C]$
denotes the subgraph of $G$ induced by $C$.
\end{proposition}

We now consider blow-ups of directed graphs.
Given a loopless directed graph $H$ and positive
integers $(n_v)_{v\in V(H)}$, a \emph{blow-up}
of $H$ is obtained by replacing each vertex
$v$ by an independent set $V_v$ of size $n_v$,
and each arc $v\to w$ by all possible arcs
from $V_v$ to $V_w$.
When $n_v=n$ for all $v\in V(H)$, we denote
the resulting graph by $H(n)$.

Let $\overrightarrow C_r$ denote the directed
cycle on $r$ vertices, with arcs $i\to i+1$,
where the indices are taken modulo $r$.
We denote its blow-up with part sizes
$n_1,\ldots,n_r$ by
$\overrightarrow C_r(n_1,\ldots,n_r)$.
For $r=2$, this is precisely the complete
bipartite graph $K_{n_1,n_2}$, with each edge
regarded as two arcs of opposite orientations.

Blow-ups of directed cycles were studied by Alon,
Ben-Eliezer, Shangguan and Tamo \cite{alon},
who showed that, for every fixed $r\geq 3$,
\[
    \HG(\overrightarrow C_r(n,\ldots,n))
    \geq n^{1/r-o(1)}.
\]
We next relate the hat-guessing number of a blow-up
of a directed cycle to that of a complete
bipartite graph.

\begin{theorem}\label{cycle}
For every $r\geq 2$ and all positive integers
$n_1,\ldots,n_r$, we have
\[
    \HG(\overrightarrow C_r(n_1,\ldots,n_r))
    \leq\min_{i\in[r]}\HG(K_{n_i,n_{i+1}}),
\]
where the indices are taken modulo $r$.
\end{theorem}

\begin{proof}
Fix $i\in[r]$, and suppose that
$q>\HG(K_{n_i,n_{i+1}})$.
Consider an arbitrary strategy on
$\overrightarrow C_r(n_1,\ldots,n_r)$
with $q$ colours.

Let $V_j$ denote the part of size $n_j$.
Notice that the guessing function of every
player in $V_i$ depends only on the colour
assignment in $V_{i+1}$.
Since $q>\HG(K_{n_i,n_{i+1}})$,
Lemma~\ref{crit} implies that there exists
a blocking set $S\subseteq Q^{V_{i+1}}$
such that none of the guessing functions
of the players in $V_i$ takes all values in
$Q$ on $S$.

For each vertex in $V_i$, choose a colour
that its guessing function does not take
on $S$, and assign a hat of this colour
to that vertex.
Thus, regardless of which assignment from
$S$ is eventually chosen for $V_{i+1}$,
all the players in $V_i$ guess incorrectly.

Next, assign hats to the vertices in $V_{i-1}$,
then to those in $V_{i-2}$, and so on,
proceeding backwards around the cycle until
only $V_{i+1}$ remains unassigned.
At each step, the hats seen by the players
in the current part have already been assigned.
We can therefore choose their hats so that
all their guesses are incorrect.

Finally, the guesses of the players in
$V_{i+1}$ are determined.
Let $z\in Q^{V_{i+1}}$ be their vector of guesses.
Since $S$ is blocking, there exists an
$x\in S\cap C_z$.
Assigning hats according to $x$ to the
vertices of $V_{i+1}$ ensures that all the
players in $V_{i+1}$ guess incorrectly.
By the choice of the hats in $V_i$, all the
players in $V_i$ also guess incorrectly.

Thus, there is no winning strategy with
$q$ colours, and hence
\[
    \HG(\overrightarrow C_r(n_1,\ldots,n_r))
    \leq\HG(K_{n_i,n_{i+1}}).
\]
Since $i\in[r]$ was arbitrary, the result
follows.
\end{proof}

Combining Theorem~\ref{cycle} with our bounds
for complete bipartite graphs gives the
following corollary.

\begin{corollary}\label{cycasymp}
For all integers $r\geq 2$, we have
\[
    \HG(\overrightarrow C_r(n,\ldots,n))
    \leq(1+o(1))\frac{n}{\log\log n}
\]
as $n\to\infty$, where the error term is
independent of $r$.

For unequal part sizes, Theorem~\ref{main}
can be applied to any pair of consecutive
parts satisfying its growth assumptions.
The same bounds hold for directed graphs
obtained by deleting arcs from the
corresponding blow-ups.
\end{corollary}

\begin{proof}
By Theorem~\ref{cycle} and
Corollary~\ref{balanced},
\[
    \HG(\overrightarrow C_r(n,\ldots,n))
    \leq\HG(K_{n,n})
    \leq(1+o(1))\frac{n}{\log\log n}.
\]
Since the bound on $\HG(K_{n,n})$ does not
depend on $r$, neither does the error term.
The assertion for unequal part sizes follows
similarly from Theorems~\ref{cycle}
and~\ref{main}.
Finally, deleting arcs cannot increase the
hat-guessing number.
\end{proof}

Notice that the bound in Corollary~\ref{cycasymp}
holds regardless of the length of the directed
cycle.
In the undirected case, a blow-up of an even
cycle is bipartite.
Thus, for any fixed even integer $r\geq 4$,
Corollary~\ref{allbip} gives
\[
    \HG(C_r(n,\ldots,n))
    \leq(1+o(1))\frac{rn}{2\log\log n},
\]
where $C_r(n,\ldots,n)$ denotes the corresponding
undirected blow-up.
However, our argument does not give a
corresponding sublinear bound for blow-ups
of undirected odd cycles.

We next give a more general sufficient
condition for obtaining sublinear upper
bounds for directed graphs.
In the proof of Theorem~\ref{cycle}, we
constructed a blocking set for one part
and then assigned hats to the vertices
in the remaining parts in a suitable order.
We extend this argument by constructing blocking
sets for several disjoint independent sets of
vertices and allowing the players outside these
sets to see additional hats.

\begin{definition}\label{decomp}
A partition
\[
    V(G)=B\sqcup A_1\sqcup\cdots\sqcup A_s,
\]
where $s\geq 1$ and each $A_j$ is nonempty,
is called \emph{admissible} if the following
conditions hold:
\begin{enumerate}
    \item Each $A_j$ is independent, and the
    directed graph on $[s]$ with an arc
    $i\to j$, for $i\ne j$, whenever there
    is an arc from $A_i$ to $A_j$ in $G$,
    is acyclic.

    \item The induced subgraph $G[B]$
    is acyclic.

    \item Every vertex in $B$ has
    out-neighbours in at most one of
    the blocks $A_1,\ldots,A_s$.
\end{enumerate}
For such a partition and a positive integer $q$,
define
\[
    r_j(q)=
    \sum_{\substack{v\in B\\
    N^+(v)\cap A_j\ne\varnothing}}
    q^{|N^+(v)\cap B|}
\]
for all $j\in[s]$.
\end{definition}

The quantity $r_j(q)$ counts the functions obtained
by considering every player in $B$ who sees at least
one hat in $A_j$ and every possible colour assignment
to the hats that this player sees in $B$.
We regard these functions as guessing functions from $Q^{A_j}$
to $Q$ by having them ignore the coordinates
corresponding to hats that the player cannot see.
Some of the resulting functions may coincide.

\begin{theorem}\label{block}
Let $G$ be a directed graph with an admissible
partition, and put
\[
    d=\max_{j\in[s]}|A_j|,
    \qquad
    r=\max\left\{1,\max_{j\in[s]}r_j(q)\right\}.
\]
Suppose that $q$ satisfies the hypotheses of
Lemma~\ref{count} and condition \eqref{test}
with these values of $d$ and $r$.
Then $\HG(G)<q$.

Alternatively, if $2\leq q\leq d$, it suffices
that condition \eqref{test0} holds with
these values of $d$ and $r$.
\end{theorem}

\begin{proof}
Fix an arbitrary strategy on $G$ with $q$
colours.
We first construct a blocking set
$S_j\subseteq Q^{A_j}$ for each $j\in[s]$.

Fix $j\in[s]$, and consider a vertex $v\in B$
such that $N^+(v)\cap A_j\ne\varnothing$.
By the third condition in
Definition~\ref{decomp}, the player at $v$
does not see any hats in the other
blocks $A_i$.

For each possible colour assignment
$b\in Q^{N^+(v)\cap B}$, define a function
\[
    f_{v,b}:Q^{A_j}\to Q
\]
as follows.
For $x\in Q^{A_j}$, let $f_{v,b}(x)$ be the
guess of the player at $v$ when the colour
assignment on $N^+(v)\cap B$ is $b$ and
that on $A_j$ is $x$.
The function ignores the coordinates of
$A_j$ corresponding to hats that the player
at $v$ cannot see.

For each such vertex $v$, there are
$q^{|N^+(v)\cap B|}$ choices of $b$.
Thus, the total number of functions obtained
in this way is at most $r_j(q)\leq r$.

Since $|A_j|\leq d$, we can apply
Theorem~\ref{family} by extending these
functions to $Q^d$, if necessary, as
described in Section~\ref{sec:sampling}.
We obtain a blocking set
$S_j\subseteq Q^{A_j}$ such that
$f_{v,b}(S_j)\ne Q$ for every function
$f_{v,b}$ considered above.

If $r_j(q)=0$, we may simply take
$S_j=Q^{A_j}$, which is blocking since
$q\geq 2$.
We choose such a blocking set separately
for each $j\in[s]$.

We now assign hats to the vertices in $B$.
Since $G[B]$ is acyclic, we may assign
these hats in reverse topological order
of $G[B]$.
When assigning a hat to a vertex $v\in B$,
all the hats that it can see in $B$ have
already been assigned.
Let $b$ denote their colour assignment.

If $v$ has an out-neighbour in some block
$A_j$, then the possible guesses of the
player at $v$, as the assignment on $A_j$
ranges over $S_j$, form the set
$f_{v,b}(S_j)$.
Since this is a proper subset of $Q$,
we can choose a colour for the hat at $v$
outside this set.

If $v$ has no out-neighbours in any of the
blocks $A_j$, its guess is already determined,
so we choose a colour different from its
guess.

In this way, we assign hats to all the
vertices in $B$ so that all the players
in $B$ guess incorrectly, regardless of
which assignments from the sets $S_j$
are subsequently chosen.

It remains to assign hats to the vertices
in the blocks $A_1,\ldots,A_s$.
By the first condition in
Definition~\ref{decomp}, the directed graph
describing the arcs between these blocks
is acyclic.
We assign hats to the blocks in reverse
topological order of this directed graph.

When considering a block $A_j$, all the
hats seen by the players in $A_j$ have
already been assigned, since $A_j$ is
independent.
Let $z_j\in Q^{A_j}$ be their vector of
guesses.
Since $S_j$ is blocking, there exists
an assignment $x_j\in S_j\cap C_{z_j}$.
Assigning hats according to $x_j$ ensures
that all the players in $A_j$ guess
incorrectly.

Repeating this procedure for all the blocks,
we obtain an assignment on which every
player in $G$ guesses incorrectly.
Hence, $\HG(G)<q$.

The proof under the alternative condition
\eqref{test0} is identical, except that
we use Proposition~\ref{orientation}
instead of Theorem~\ref{family} to
construct the blocking sets $S_j$.
\end{proof}

Notice that the quantity $r$ in
Theorem~\ref{block} depends on $q$.
To obtain an asymptotic bound, we replace
it by an upper bound that is independent
of $q$.

\begin{corollary}\label{blockasymp}
Consider a sequence of directed graphs $G$
with admissible partitions, where
$d=\max_{j\in[s]}|A_j|\to\infty$.
Put
\[
    R_d=\max\left\{
    d,\,
    \max_{j\in[s]}
    \sum_{\substack{v\in B\\
    N^+(v)\cap A_j\ne\varnothing}}
    d^{|N^+(v)\cap B|}
    \right\}.
\]
If $\log R_d=o(d)$, then
\[
    \HG(G)\leq(1+o(1))
    \frac{d}{\log\log(d/\log R_d)}
    =o(d).
\]
In particular, if
$\log R_d=O(d^{1-\delta})$ for some
fixed $\delta>0$, then
\[
    \HG(G)\leq(1+o(1))\frac{d}{\log\log d}.
\]
\end{corollary}

\begin{proof}
For all $q\leq d$ and all $j\in[s]$, we have
\[
    r_j(q)\leq
    \sum_{\substack{v\in B\\
    N^+(v)\cap A_j\ne\varnothing}}
    d^{|N^+(v)\cap B|}
    \leq R_d.
\]
Consequently,
\[
    \max\left\{1,\max_{j\in[s]}r_j(q)\right\}
    \leq R_d.
\]

Since $d\leq R_d$ and $\log R_d=o(d)$,
we may choose $q$ as in the proof of
Theorem~\ref{main}, with $m=d$ and $n=R_d$.
For this choice of $q$, condition
\eqref{test} holds with $r=R_d$,
and hence also with
$r=\max\{1,\max_{j\in[s]}r_j(q)\}$.
Thus, Theorem~\ref{block} gives
\[
    \HG(G)\leq(1+o(1))
    \frac{d}{\log\log(d/\log R_d)}.
\]
Since $\log R_d=o(d)$, the denominator
tends to infinity, proving that
$\HG(G)=o(d)$.

If $\log R_d=O(d^{1-\delta})$, then,
as in the proof of Corollary~\ref{balanced},
we have
\[
    \log\log(d/\log R_d)
    =\log\log d+O(1).
\]
This proves the last assertion.
\end{proof}

For example, suppose that for each block
$A_j$, at most $p$ vertices in $B$ have
out-neighbours in $A_j$, and that each of
these vertices has at most $h$ out-neighbours
in $B$.
Then
\[
    R_d\leq\max\{d,pd^h\}.
\]
Consequently, if
$\log p+h\log d=o(d)$, where $p\geq 1$,
Corollary~\ref{blockasymp} gives a sublinear
bound on $\HG(G)$.

We next apply Corollary~\ref{blockasymp}
to blow-ups of certain fixed directed graphs.

\begin{corollary}\label{base}
Let $H$ be a fixed loopless directed graph,
and suppose that there exists a nonempty
set $I\subseteq V(H)$ such that both
$H[I]$ and $H[V(H)\setminus I]$ are acyclic.
Suppose further that every vertex outside
$I$ that has an out-neighbour in $I$ has
exactly one out-neighbour in $I$ and no
out-neighbours outside $I$.
Then
\[
    \HG(H(n))\leq
    (1+o(1))\frac{n}{\log\log n}
\]
as $n\to\infty$.
\end{corollary}

\begin{proof}
For every $v\in I$, take the part of $H(n)$
corresponding to $v$ as one of the blocks
$A_1,\ldots,A_s$.
Let $B$ be the union of the remaining parts.

Since $H[I]$ is acyclic, the directed graph
describing the arcs between the blocks
$A_j$ is acyclic.
Moreover, since a blow-up of an acyclic
directed graph is acyclic, the induced
subgraph on $B$ is acyclic.

By the assumption on the out-neighbours
of the vertices outside $I$, every vertex
in $B$ has out-neighbours in at most one
of the blocks $A_j$.
Thus, the partition is admissible.

Furthermore, every vertex in $B$ that has
an out-neighbour in some block $A_j$ has
no out-neighbours in $B$.
Since each block has size $n$, we have
$d=n$ and
\[
    R_n\leq\max\{n,|V(H)|n\}.
\]
As $H$ is fixed, it follows that
$\log R_n=O(\log n)$.
The result now follows from
Corollary~\ref{blockasymp}.
\end{proof}

For example, consider the directed graph
$H$ on four vertices $a,b,c,d$, with arcs
\[
    a\to b,\quad a\to c,\quad
    b\to d,\quad c\to d,\quad d\to a.
\]
Taking $I=\{a\}$, we see that
Corollary~\ref{base} applies, so
\[
    \HG(H(n))\leq
    (1+o(1))\frac{n}{\log\log n}.
\]
Thus, the result applies not only to blow-ups
of directed cycles, but also to blow-ups
of certain other strongly connected
directed graphs.

Suppose instead that $I$ is a nonempty
independent subset of $V(H)$, that
$H[V(H)\setminus I]$ is acyclic, and that
every vertex outside $I$ that has an
out-neighbour in $I$ has no out-neighbours
outside $I$.
In this case, we may allow such a vertex
to have more than one out-neighbour in $I$.

Indeed, we can take the union of all the
parts corresponding to vertices of $I$
as a single independent block.
This gives an admissible partition with
$d=|I|n$.
As above, every vertex in $B$ that has
an out-neighbour in this block has no
out-neighbours in $B$, so
$R_d=O(n)$.
Since $H$ is fixed, Corollary~\ref{blockasymp}
gives
\[
    \HG(H(n))\leq
    (1+o(1))\frac{|I|n}{\log\log n}.
\]

The results above also apply to blow-ups
with unequal part sizes whenever the
corresponding assumptions of
Theorem~\ref{block} or
Corollary~\ref{blockasymp} are satisfied.

Finally, Proposition~\ref{scc} allows us
to apply these results separately to the
strongly connected components of a
directed graph.
For example, every strongly connected
component of a loopless directed graph
of maximum out-degree at most one is
either a directed cycle or an isolated
vertex within that component.
Consequently, the bounds for blow-ups
of directed cycles also give bounds for
blow-ups of such directed graphs.

\section{A limitation of the method and related questions}
\label{sec:barrier}

In this section, we show that the asymptotic
bound in Theorem~\ref{main} cannot be improved
within the probabilistic framework used in
its proof.

Recall that in the proof of
Theorem~\ref{family}, we constructed a
random set $S$ whose distribution was
independent of the functions
$f_1,\ldots,f_r:Q^d\to Q$.
We then showed that, with positive
probability, $S$ was blocking and none
of the functions $f_i$ took all values
in $Q$ on $S$.

More generally, given $d,q,r$, suppose
that there exists a probability distribution
$\mathcal D$ on the blocking subsets
of $Q^d$ such that
\begin{equation}\label{univgoal}
    \sup_{f:Q^d\to Q}
    \Pp_{S\sim\mathcal D}(f(S)=Q)
    <\frac1r.
\end{equation}
Then, for any family of $r$ functions
$f_1,\ldots,f_r:Q^d\to Q$, the union
bound gives
\[
    \Pp_{S\sim\mathcal D}
    (f_i(S)=Q\text{ for some }i\in[r])<1.
\]
Hence, there exists a blocking set
$S\subseteq Q^d$ such that
$f_i(S)\ne Q$ for all $i\in[r]$.

In fact, under condition \eqref{test},
the construction used in the proof
of Theorem~\ref{family}, conditioned
on $S$ being blocking, satisfies
\eqref{univgoal}.
Indeed, by the proof of that theorem,
we have
\[
    \Pp(S\text{ is blocking})\geq\frac34,
\]
and, for every fixed function $f:Q^d\to Q$,
\[
    \Pp(f(S)=Q)\leq\frac1{2r}.
\]
Consequently,
\begin{align*}
    \Pp(f(S)=Q\mid S\text{ is blocking})
    &\leq
    \frac{\Pp(f(S)=Q)}
    {\Pp(S\text{ is blocking})}\\
    &\leq\frac{2}{3r}<\frac1r.
\end{align*}
Thus, it is natural to ask for which values
of $d,q,r$ a distribution satisfying
\eqref{univgoal} exists.

A blocking set in $Q^d$ is precisely a total
dominating set in the direct product of $d$
copies of $K_q$, that is, the graph on $Q^d$
in which two vertices are adjacent if they
differ in every coordinate.
Lemmas~3.2 and~3.3 of \cite{alon} imply that
every blocking set has size at least $d+1$
and greater than $e^{d/q}$.
Related lower bounds and coordinate-by-coordinate
arguments appear in \cite{defant,vemuri}.
We use a similar argument to establish the
following bound.

\begin{lemma}\label{size}
For $2\leq q\leq d$, every blocking set
$S\subseteq Q^d$ satisfies
\[
    |S|>(q-1)
    \left(\frac{q}{q-1}\right)^{d-q+1}.
\]
In particular,
\[
    \frac{|S|}{q}>\frac{e^{d/q}}{2e}.
\]
\end{lemma}

\begin{proof}
Suppose, for a contradiction, that
\[
    |S|\leq(q-1)
    \left(\frac{q}{q-1}\right)^{d-q+1}.
\]
We will construct a cube $C_z$ such
that $S\cap C_z=\varnothing$.

Choose $z_1$ to be a most frequent
colour in the first coordinate among
the elements of $S$.
Since there are $q$ colours, at least
$|S|/q$ elements of $S$ have this
colour in their first coordinate.
None of these elements can belong
to $C_z$.

Next, choose $z_2$ to be a most
frequent colour in the second
coordinate among the elements of
$S$ that have not yet been excluded.
This excludes at least a $1/q$
fraction of the remaining elements.

Continue in this way for the first
$d-q+1$ coordinates.
After these coordinates have been
chosen, the number of elements of
$S$ that have not yet been excluded
is at most
\[
    |S|\left(1-\frac1q\right)^{d-q+1}
    \leq q-1.
\]
There are $q-1$ coordinates of $z$
that have not yet been chosen.
We assign a different unused
coordinate to each of the remaining
elements of $S$, and choose the
colour of $z$ in that coordinate
to agree with the corresponding
element of $S$.
Any remaining coordinates of $z$
can be chosen arbitrarily.

The resulting cube $C_z$ contains
no element of $S$, contradicting
the assumption that $S$ is blocking.
This proves the first inequality.

For the second inequality, observe that
\[
    -\log\left(1-\frac1q\right)\geq\frac1q.
\]
Hence,
\begin{align*}
    \frac{|S|}{q}
    &>\frac{q-1}{q}
    \left(1-\frac1q\right)^{-(d-q+1)}\\
    &\geq\frac{q-1}{q}
    \exp\left(\frac{d-q+1}{q}\right)\\
    &\geq\frac12 e^{d/q-1}\\
    &=\frac{e^{d/q}}{2e},
\end{align*}
where we used $q\geq 2$.
\end{proof}

The idea of the argument that follows
is that, if $q$ is sufficiently small
relative to $d$, then every blocking
set must be sufficiently large for a
random function to take all values in
$Q$ on that set with reasonably large
probability.

Let $F:Q^d\to Q$ be a random function
obtained by choosing $F(x)$ independently
and uniformly from $Q$ for every
$x\in Q^d$.
For any fixed set $S\subseteq Q^d$,
the probability that a particular
colour does not belong to $F(S)$ is
$(1-1/q)^{|S|}$.
Thus, by the union bound,
\begin{equation}\label{missingcolour}
    \Pp_F(F(S)\ne Q)
    \leq q\left(1-\frac1q\right)^{|S|}
    \leq qe^{-|S|/q}.
\end{equation}

As we will see later, for balanced or moderately
unbalanced complete bipartite graphs, this estimate,
together with Lemma~\ref{size}, is sufficient to
establish the desired limitation of our method.
For the full range of part sizes in
Theorem~\ref{main}, however, we need a lower bound
on the probability that $F(S)=Q$, even when this
probability tends to zero.
To this end, we use the standard technique of
Poissonisation (see, for example,
\cite[Section~5.4]{mitzenmacher}) to prove
the following lemma.

\begin{lemma}\label{poisson}
Let $p_q(k)$ denote the probability
that $k$ independent and uniformly chosen colours
from a set $Q$ of $q\geq 2$ colours
include every colour in $Q$.
Then, for sufficiently large $k/q$,
\[
    p_q(k)\geq
    \frac12\exp\left(-2qe^{-k/(2q)}\right).
\]
\end{lemma}

\begin{proof}
Consider an infinite sequence of colours chosen
independently and uniformly from $Q$, and let $L$
be a random variable, independent of these choices,
having the Poisson distribution with mean $k/2$. Thus,
\[
    \Pp(L=l)=e^{-k/2}\frac{(k/2)^l}{l!}
\]
for all nonnegative integers $l$.

For each $a\in Q$, let $C_a$ denote the number
of occurrences of colour $a$ among the first
$L$ colours in the sequence.

We claim that the random variables
$C_a$, for $a\in Q$, are mutually
independent and each has the Poisson
distribution with mean $k/(2q)$.
Indeed, for any choice of nonnegative
integers $c_a$, for $a\in Q$, put
$l=\sum_{a\in Q}c_a$.
Then
\begin{align*}
    &\Pp(C_a=c_a\text{ for all }a\in Q)\\
    &\qquad=
    \Pp(L=l)\frac{l!}{\prod_{a\in Q}c_a!}
    q^{-l}\\
    &\qquad=
    e^{-k/2}\frac{(k/2)^l}{l!}
    \frac{l!}{\prod_{a\in Q}c_a!}q^{-l}\\
    &\qquad=
    \prod_{a\in Q}
    \left(e^{-k/(2q)}
    \frac{(k/(2q))^{c_a}}{c_a!}\right).
\end{align*}
The last expression is the product
of the probabilities that independent
Poisson random variables with mean
$k/(2q)$ take the corresponding values.
This proves the claim.

In particular, for each $a\in Q$,
\[
    \Pp(C_a=0)=e^{-k/(2q)}.
\]
By independence, the probability that every colour
occurs among the first $L$ colours in the sequence is
\[
    \Pp(C_a\geq 1\text{ for all }a\in Q)
    =\left(1-e^{-k/(2q)}\right)^q.
\]

If every colour occurs among the first $L$ colours
in the sequence and $L\leq k$, then every colour
also occurs among the first $k$ colours.
Consequently,
\begin{equation}\label{poissoncompare}
    p_q(k)\geq
    \left(1-e^{-k/(2q)}\right)^q
    -\Pp(L>k).
\end{equation}

We next estimate $\Pp(L>k)$.
By the definition of the Poisson
distribution,
\begin{align*}
    \Ee[2^L]
    &=\sum_{l=0}^{\infty}
    2^l e^{-k/2}\frac{(k/2)^l}{l!}\\
    &=e^{-k/2}\sum_{l=0}^{\infty}
    \frac{k^l}{l!}\\
    &=e^{k/2}.
\end{align*}
Thus, by Markov's inequality,
\[
    \Pp(L>k)
    \leq\frac{\Ee[2^L]}{2^k}
    =\exp\left(-k\left(\log 2-\frac12\right)\right).
\]

On the other hand, using the inequality
$\log(1-u)\geq-2u$ for $0\leq u\leq 1/2$,
we obtain, for sufficiently large $k/q$,
\[
    \left(1-e^{-k/(2q)}\right)^q
    \geq\exp\left(-2qe^{-k/(2q)}\right).
\]
Moreover,
\begin{align*}
    \frac{\exp(-k(\log 2-1/2))}
    {\exp(-2qe^{-k/(2q)})}
    &=\exp\left(-k(\log 2-1/2)+2qe^{-k/(2q)}\right)\\
    &=\exp\left(-k\left(\log 2-\frac12
    -\frac{2q}{k}e^{-k/(2q)}\right)\right)
    =o(1)
\end{align*}
as $k/q\to\infty$, since
$\log 2-1/2>0$.
In particular, for sufficiently
large $k/q$, we have
\[
    \Pp(L>k)
    \leq\frac12\exp\left(-2qe^{-k/(2q)}\right).
\]
Substituting these estimates into
\eqref{poissoncompare}, we obtain
\[
    p_q(k)\geq
    \frac12\exp\left(-2qe^{-k/(2q)}\right),
\]
as required.
\end{proof}

We are now ready to establish the
limitation of the probabilistic
framework described in \eqref{univgoal}.

\begin{theorem}\label{barrier}
Let $d\to\infty$, with $r\geq d$
and $\log r=o(d)$.
Fix $0<\delta<1$, and suppose that
\[
    2\leq q\leq
    \frac{(1-\delta)d}{\log\log(d/\log r)}.
\]
Then, for every probability
distribution $\mathcal D$ on the
blocking subsets of $Q^d$, we have
\[
    \sup_{f:Q^d\to Q}
    \Pp_{S\sim\mathcal D}(f(S)=Q)
    \geq\exp(-o(\log r))>\frac1r
\]
for sufficiently large $d$.
The $o(\log r)$ term is uniform in $\mathcal D$
and in $q$ in the indicated range.
In particular, \eqref{univgoal}
cannot hold.
\end{theorem}

\begin{proof}
Put
\[
    R=\log(d/\log r).
\]
Since $\log r=o(d)$, we have
$R\to\infty$.
By the assumption on $q$,
\[
    \frac{d}{q}\geq
    \frac{\log R}{1-\delta}.
\]
Let $S\subseteq Q^d$ be an arbitrary
blocking set, and put $k=|S|$.
By Lemma~\ref{size},
\begin{equation}\label{barriersize}
    \frac{k}{q}>
    \frac{e^{d/q}}{2e}
    \geq\frac{R^{1/(1-\delta)}}{2e}.
\end{equation}
In particular, $k/q\to\infty$,
uniformly over all blocking sets $S$
and all $q$ in the indicated range.

Choose a random function $F:Q^d\to Q$
by choosing $F(x)$ independently
and uniformly from $Q$ for every
$x\in Q^d$.
Since $S$ is fixed and has size $k$,
we have
\[
    \Pp_F(F(S)=Q)=p_q(k).
\]
Consequently, by Lemma~\ref{poisson},
\[
    \Pp_F(F(S)=Q)
    \geq\frac12\exp\left(-2qe^{-k/(2q)}\right).
\]

We claim that
\[
    qe^{-k/(2q)}=o(\log r).
\]
Indeed, for sufficiently large $d$, we have
$q\leq d$, so it follows from \eqref{barriersize}
that
\begin{align*}
    \frac{qe^{-k/(2q)}}{\log r}
    &\leq\frac{d}{\log r}
    \exp\left(
    -\frac{R^{1/(1-\delta)}}{4e}
    \right)\\
    &=\exp\left(
    R-\frac{R^{1/(1-\delta)}}{4e}
    \right)\\
    &=o(1),
\end{align*}
where the last equality follows from
$R\to\infty$ and $1/(1-\delta)>1$.
This proves the claim.

Since also $\log r\to\infty$, we obtain
\[
    \Pp_F(F(S)=Q)\geq\exp(-o(\log r)).
\]
Importantly, this lower bound holds
uniformly over all blocking sets $S$.

Now let $S$ be chosen according to
an arbitrary probability distribution
$\mathcal D$ on the blocking subsets
of $Q^d$, independently of $F$.
Averaging the preceding inequality
over $S$, we obtain
\[
    \Pp_{F,S}(F(S)=Q)
    \geq\exp(-o(\log r)).
\]
On the other hand, by averaging
first over $S$ and then over $F$,
we have
\[
    \Pp_{F,S}(F(S)=Q)
    =\Ee_F\left[
    \Pp_{S\sim\mathcal D}(F(S)=Q)
    \right].
\]
Hence, there exists a fixed function
$f:Q^d\to Q$ such that
\[
    \Pp_{S\sim\mathcal D}(f(S)=Q)
    \geq\exp(-o(\log r)).
\]
Finally, since $\log r\to\infty$,
we have
$\exp(-o(\log r))>1/r$ for sufficiently
large $d$.
This completes the proof.
\end{proof}

When $r$ is not too large relative
to $d$, we can obtain a stronger
conclusion using the simpler
estimate \eqref{missingcolour}.

\begin{corollary}\label{strongbarrier}
Under the hypotheses of
Theorem~\ref{barrier}, suppose
further that
$\log r=O(d^{1-\varepsilon})$
for some fixed $\varepsilon>0$.
Then
\[
    \sup_{f:Q^d\to Q}
    \Pp_{S\sim\mathcal D}(f(S)=Q)
    =1-o(1)
\]
for every probability distribution
$\mathcal D$ on the blocking subsets
of $Q^d$, uniformly in $\mathcal D$.
In particular, this holds when $r=d$.
\end{corollary}

\begin{proof}
As in the proof of
Theorem~\ref{barrier}, put
$R=\log(d/\log r)$.
By the additional assumption on $r$,
we have
\[
    R\geq\varepsilon\log d-O(1).
\]
For an arbitrary blocking set
$S\subseteq Q^d$ and a random
function $F:Q^d\to Q$ as above,
it follows from \eqref{missingcolour}
and \eqref{barriersize} that, for sufficiently
large $d$,
\begin{align*}
    \Pp_F(F(S)\ne Q)
    &\leq qe^{-|S|/q}\\
    &\leq d\exp\left(
    -\frac{R^{1/(1-\delta)}}{2e}
    \right)\\
    &=o(1).
\end{align*}
Indeed, since $1/(1-\delta)>1$,
the quantity $R^{1/(1-\delta)}$
grows faster than $\log d$.
The error term is independent of
the choice of the blocking set $S$.

Consequently,
\[
    \Pp_F(F(S)=Q)\geq 1-o(1).
\]
Averaging over $S\sim\mathcal D$
and then over $F$, as in the proof
of Theorem~\ref{barrier}, we obtain
a fixed function $f:Q^d\to Q$
such that
\[
    \Pp_{S\sim\mathcal D}(f(S)=Q)
    \geq 1-o(1).
\]
This proves the result.
\end{proof}

Taking $d=m$ and $r=n$,
Theorem~\ref{barrier} shows that
the probabilistic framework
described in \eqref{univgoal}
cannot give an asymptotic upper
bound of the form
\[
    \HG(K_{m,n})\leq
    \frac{(1-\delta)m}
    {\log\log(m/\log n)}
\]
for any fixed $\delta>0$ throughout
the range of part sizes considered
in Theorem~\ref{main}.

One may ask whether a better bound
can be obtained by constructing
blocking sets on the larger part
instead.
We next show that the same
obstruction holds in that case.

Suppose that $m\leq n$,
$\log n=o(m)$, and
\[
    2\leq q\leq
    \frac{(1-\delta)m}
    {\log\log(m/\log n)}
\]
for some fixed $0<\delta<1$.
Let $S\subseteq Q^n$ be an arbitrary
blocking set.
By Lemma~\ref{size},
\[
    \frac{|S|}{q}>
    \frac{e^{n/q}}{2e}.
\]

Suppose first that $n\leq m^2$.
Since $m\leq n\leq m^2$, we have
\[
    \log(m/\log n)\sim\log m.
\]
Moreover, since $n\geq m$,
\begin{align*}
    \frac{|S|}{q}
    &>\frac{e^{n/q}}{2e}\\
    &\geq\frac{e^{m/q}}{2e}\\
    &\geq
    \frac{(\log(m/\log n))^{1/(1-\delta)}}{2e}.
\end{align*}
It follows that there exists a
constant $c>0$ such that
\[
    \frac{|S|}{q}
    \geq c(\log m)^{1/(1-\delta)}
\]
for sufficiently large $m$.

Suppose now that $n>m^2$.
Since $q\leq m$ for sufficiently
large $m$, we have $n/q>m$.
Hence,
\[
    \frac{|S|}{q}>\frac{e^m}{2e}.
\]
In either case, since $q\leq m$,
we obtain
\[
    qe^{-|S|/q}=o(1).
\]

Thus, by \eqref{missingcolour},
for a random function $F:Q^n\to Q$,
we have $\Pp_F(F(S)=Q)\geq 1-o(1)$,
uniformly over all blocking sets
$S\subseteq Q^n$.
Averaging as before, we conclude
that for every probability
distribution $\mathcal D$ on the
blocking subsets of $Q^n$, there
exists a function $f:Q^n\to Q$
such that
\[
    \Pp_{S\sim\mathcal D}(f(S)=Q)
    \geq 1-o(1).
\]
In particular, the analogue of
\eqref{univgoal} with $d=n$ and
$r=m$ cannot hold.

Therefore, constructing blocking
sets on the larger part instead
of the smaller part does not
allow us to improve the leading
constant in our asymptotic bound
within this framework.

It is worth emphasising that these
results do not give a lower bound
on the hat-guessing number of
$K_{m,n}$ itself.
Our result only shows a limitation of the method
in which the distribution of blocking sets is
chosen independently of the guessing functions.
More specifically, it shows that the union bound
based on the individual surjectivity probabilities
cannot give a better asymptotic bound.
It does not rule out choosing a clever
blocking set that depends on the
particular guessing functions,
or exploiting the joint behaviour
of these functions.

The gap between the lower bound
of Alon, Ben-Eliezer, Shangguan
and Tamo \cite{alon} and the
upper bound in Theorem~\ref{main}
remains large.
In particular, the following
problem seems natural.

\begin{problem}\label{polyproblem}
Does there exist an absolute constant
$\varepsilon>0$ such that
\[
    \HG(K_{n,n})=O(n^{1-\varepsilon})?
\]
\end{problem}

By Theorem~\ref{barrier}, a positive
answer to this problem cannot be
obtained merely by improving the
estimates in the probabilistic
argument described in
\eqref{univgoal}.

Another natural question concerns blow-ups of
undirected odd cycles. In contrast to the directed
case, each player in such a blow-up sees hats in
two different parts, so the argument used in the
proof of Theorem~\ref{cycle} does not apply.
The simplest instance is the complete tripartite
graph $K_{n,n,n}$.
It is known that
$\HG(K_{n,n,n})\geq n^{2/3-o(1)}$
(see \cite[Theorem~1.2]{alon}).
Can we show that $\HG(K_{n,n,n})=o(n)$?

\paragraph{Declaration of AI use.}
OpenAI's ChatGPT was used to assist in the development
of the mathematical arguments, the literature search
and the preparation of the manuscript.
The author takes full responsibility for the
mathematical content and the final version of the paper.

\end{document}